\documentclass[11pt]{iopart}
\usepackage{amsmath}
\usepackage{fullpage}
\usepackage{graphicx}
\usepackage{amssymb,graphicx,float,epsf}
\usepackage{float}
\restylefloat{figure}

\def\NN{\mathbb N}
\def\RR{\mathbb R}
\newcommand{\qed}{\hbox to 0pt{}\hfill$\rlap{$\sqcap$}\sqcup$\vspace{2mm}}
\newtheorem{theorem}{Theorem}[section]

\numberwithin{equation}{section}

\newtheorem{uess}[theorem]{Lemma}
\newtheorem{guess}[theorem]{Theorem}
\newtheorem{prop}[theorem]{Proposition}
\newtheorem{remark}[theorem]{Remark}
\newtheorem{corol}[theorem]{Corollary}
\newtheorem{example}[theorem]{Example}
\newtheorem{definition}[theorem]{Definition}

\begin{document}

\title[Global attractivity of non-autonomous neural networks]{On the global attractivity of non-autonomous neural networks 
with a distributed delay}

\author{Leonid Berezansky$^1$ and Elena Braverman$^2$}
\address{$^1$ Dept. of Math, Ben-Gurion University of Negev, Beer-Sheva 84105, 
Israel}
\ead{brznsky@math.bgu.ac.il}
\address{$^2$ Dept. of Math \& Stats, University of Calgary,
2500 University Dr. NW, Calgary, AB, Canada T2N 1N4}
\ead{maelena@math.ucalgary.ca}



\begin{abstract}
We consider a system of several nonlinear equations with a distributed delay and obtain absolute asymptotic stability conditions, independent of the delay.
The ideas of the proofs are based on the notion of a strong attractor. The results are applied to Hopfield neural networks,  Nicholson's blowflies type system, and compartment models of population dynamics.
\end{abstract}



\maketitle

\section{Introduction}

In many applications, systems experiencing delay in fact involve distributed delays. If the rate of change depends on the past, for example, reporting delays in economics or maturation delays in ecological systems, assuming constant concentrated delays is a significant simplification.
Generally, the dependency will be on a certain segment of prehistory of the process, and, in certain cases, on the whole of it, leading to infinite delays. This is the reason why stability of systems with a distributed delay has been intensively investigated, let us mention some recent publications 
\cite{Aliseyko,Automatica2012,Oliveira,Faria_JDE_2017,Glizer,Beretta,Liu,Muroya,Solomon,Xu}, see also references therein.
In biological models \cite{Gourley,Hattaf,Yang}, a distributed delay sometimes is accepted by default.
A neural network \cite{Hopfield} is one of the most important applications where delays, in particular distributed, occur.

For scalar equations with a distributed delay, stability, either dependent or independent of the delay distribution, has been studied by several authors \cite{Nonlin2013,JMAA2014,ZAA2019,Campbell,Yuan_2011}.
A particular case of a system of two equations was explored in \cite{Nonlin2015,Cheng}.

The purpose of the present paper is to obtain delay-independent stability conditions for a system 
of differential equations with a distributed delay
\begin{equation} 
\label{1a}
\frac{dX}{dt} = G(t) \left[ \int_{H(t)}  d_{\tau} R(t,\tau)  F(X(\tau)) - X(t) \right],
\end{equation} 
where $X: \RR \to \RR^s$, $s \in {\mathbb N}$, $X$ is a column vector function $X=(x_1,x_2,\dots x_s)$, the column vector  function of $s$ variables 
$F:\RR^s \to \RR^{s}$ is
\begin{equation}
\label{F_def}
F=(f_1,f_2, \dots f_s),
\end{equation}
$G: \RR \to \RR_+^s \times \RR_+^s$, $\RR^+=[0,+\infty)$, and
$R:\RR^{s+1} \to \RR^s \times \RR^s$ are matrix functions: $G$ is diagonal
with $g_i$ on the diagonal, and $r_{ij}(t,\tau)$ are entries of $R$. We consider the case when the $s$-dimensional domain 
$H(t)$ is $H(t)=[h_1(t),t] \times [h_2(t),t] \times \dots \times [h_s(t),t]$,
$h_j(t)\leq t$, and the volume integral can be iterated, where
system \eqref{1a} can be rewritten as a collection of $s$
equations for $i=1,2, \dots, s$,
\begin{equation}
\label{intro1}
\frac{dx_i}{dt} = g_i(t) \left[ \int_{h_{i1}(t)}^t \! \!\!\! \!\!d_{\tau_1} r_{i1}(t,\tau_1) 
\dots \int_{h_{is}(t)}^t \!\!\! \!\!\! d_{\tau_s} r_{is}(t,\tau_s) f_i(x_1(\tau_1),x_2(\tau_2), \dots, x_s(\tau_s))  - x_i(t) \right] .
\end{equation} 
In particular, let $h_{ij}(t) \leq t$, $i,j=1, \dots, s$ be measurable functions, and $r_{ij}$ be step functions taking the value of one on half-open intervals $(h_{ij}(t),+\infty)$: 
\begin{equation*}
r_{ij} (t,\zeta) = \chi_{(h_{ij}(t),\infty)} (\zeta), ~~~ 
\chi_J (\zeta) := \left\{ \begin{array}{ll} 1, & \zeta \in J, \\
0, & \zeta \not\in J. \end{array} \right.
\end{equation*} 
Then, \eqref{intro1} has the form
\begin{equation}
\label{concentr_delay}
\frac{dx_i}{dt} = g_i(t) \left[ f_i(x_1(h_{i1}(t)),x_2(h_{i2}(t)), \dots, x_s(h_{is}(t)))  - x_i(t) \right] , ~~i=1,2, \dots, s.
\end{equation} 

The  Hopfield neural network \cite{Hopfield} 
\begin{equation} 
\label{intro2a}
x_i^{\prime}(t) = -b_i x_i(t) + \sum_{j=1}^s c_{ij} \tilde{f}_j(x_j(t-\tau_{ij})), ~~t\geq 0,~i=1, \dots, s
\end{equation}
is a particular case of \eqref{concentr_delay} for
$$
g_i(t) \equiv b_i, ~~~ f_i(x_1,x_2, \dots, x_n) = \frac{1}{b_i} \sum_{j=1}^s c_{ij} \tilde{f}_j(x_j),
~~h_{ij}(t) \equiv t- \tau_{ij}, ~~i,j=1, \dots, s.
$$
In the case of absolutely continuous in $\zeta$ functions $r_{ij} (t,\zeta)$,
\begin{equation*}
\frac{\partial }{ \partial \zeta} r_{ij} (t,\zeta) = k_{ij} (t,\zeta),~~i,j=1, \dots, s, ~~~
K_i(t,\tau_1,\tau_2, \dots, \tau_s) = \prod_{j=1}^s k_{ij} (t,\tau_j),
\end{equation*}
\eqref{intro1} becomes
\begin{equation}
\label{all_distr_delay}
\frac{dx_i}{dt} = g_i(t) \left[ \int\limits_{H(t)} K_i(t,\tau_1, \dots, \tau_s) f_i(x_1(\tau_1), \dots, x_s(\tau_s))~d\tau_1 \dots d\tau_s  - x_i(t) \right], ~i=1, \dots, s.
\end{equation}

In future, we consider each of the equations separately, as in \eqref{intro1}. However, due to the length of  \eqref{intro1},
a shorter notation will be used, aligned with \eqref{1a}
\begin{equation} 
\label{intro1a}
\frac{dx_i}{dt} = g_i(t) \left[ \int_{H(t)}  d_{\tau} R_i(t,\tau)  f_i(X(\tau)) - x_i(t) \right],
\end{equation} 
where $\tau=(\tau_1, \dots, \tau_s)$.

The purpose of the present paper is to explore global asymptotic stability of 
cooperative systems with distributed delays, which include \eqref{intro2a} and \eqref{all_distr_delay} as special cases.
In addition to being distributed, the delays can change with time.
Distributed delays describe a feasible fact that any interval for delay values has some probability,
such models include equations with concentrated (either constant or variable) delays. 

Compared to most previous works, main differences are outlined below.
\begin{itemize}
\item
Distributed delays can, as particular cases, include
systems with variable concentrated delays, integral terms (used in most papers on distributed delays), 
their combinations, and some other models (for example, the Cantor function as a distribution).
Moreover, argument deviations can be Lebesgue measurable, not necessarily continuous, functions.
This is the reason why methods for continuous delays do not work in this setting.
\item
Delay distributions can be non-autonomous. 
If we describe these distributions as a probability that
a delay takes a greater than a given value, this corresponds to time-dependent delay. In applications,
this allows to consider, for example, seasonal changes in delay distributions. To some extent, we explore
the most general system with a unique positive equilibrium, and justify global stability of this equilibrium,
once delays are involved only in those terms which describe cross-influences. 
The present paper answers the question when delays do not have any destabilizing effect on a non-autonomous
system.
\item
On the other hand, many of the previous papers on distributed delay describe much more complicated dynamics than
absolute global stability established in the present paper. For example, delay dependence of stability properties  
was studied in \cite{BrZhuk}, while possible multistability considered in \cite{JMAA2014}.
However, the study of systems which can be destabilized by large enough delay are not in the framework of the present paper.
Here we restrict ourselves to ``absolutely stable'' systems, where no type or size of a finite delay can destabilize it, as long as the initial conditions belong to the ``attraction set''.
\end{itemize}


The plan of the paper is as follows. After some preliminaries and an auxiliary statement in Section~\ref{sec:prelim}, we get 
stability results for systems with a distributed delay in Section~\ref{sec:stability}. 
These theorems are later applied to particular cases of neural networks and models of population dynamics in Section~\ref{sec:applications}. 
Finally, the results are discussed, and some open problems and directions of research are outlined in Section~\ref{sec:discussion}.

\section{Preliminaries}
\label{sec:prelim}

Consider a system with distributed delays \eqref{1a},
under the initial condition  
\begin{equation}
\label{2star}
X(t)=\Phi(t),~t \leq t_0,
\end{equation}
where $\Phi(t)$ is a bounded vector function.

\begin{definition}
A vector function $X(t)$ is {\bf a solution  of system} \eqref{1a},\eqref{2star} if
it satisfies \eqref{1a} for almost all $t \geq t_0$ and \eqref{2star} for $t \leq t_0$.
\end{definition}

In particular, \eqref{intro1} can be written in matrix form \eqref{1a}, where 
\begin{equation}
\label{1b}
H(t)= 
(h_1(t),t] \times (h_2(t),t] \times \cdots \times (h_s(t),t].
\end{equation}

Problems \eqref{intro1},\eqref{2star} and \eqref{1a},\eqref{2star} will be investigated under some of the following assumptions.
\begin{description}
\item{{\bf (a1)}} 
There is a domain $D \subset \RR^s$ such that all 
$f_i: D \to \RR$  are continuous functions,  $i=1, \dots, s$.

\item{{\bf (a2)}} 
Any scalar delay function $h:\RR^+ \rightarrow \RR$ considered in the paper ($h_{ij}$ in particular)
 is Lebesgue measurable, $h(t)\leq t$ and $\lim\limits_{t\rightarrow +\infty}h(t)=+\infty$.

\item{{\bf (a3)}} 
The entries of the matrix $\displaystyle R(t,\tau)=\left(  r_{ij} (t,\tau) \right)_{i,j=1}^s$,
$r_{ij}(t, \cdot)$, $i,j=1,\dots, s$   are left continuous non-decreasing functions
for any $t$, $r_{ij}(\cdot,\tau)$ are locally integrable for
any $s$, $r_{ij}(t,\tau)=0$, $\tau \leq h_i(t)$, $r_{ij}(t,t^+)=1$, $i,j=1,\dots, s$,
here all the integrals are understood in the sense
$$
\int_{h(t)}^t f(\zeta)\, dr(\zeta) = \int_{h(t)}^{t+} f(\zeta)\, dr(\zeta), ~~
\int_{h(t)}^t f(\zeta)\, d \chi_{[t,+\infty)}(\zeta) = f(t),
$$
where $u(t^+)$ is the right-side limit of the function $u$ at point $t$.

\item{{\bf (a4)}}
$G(t)$=diag$\{ g_1(t), \dots, g_s(t) \}$, $g_i(t)$ are Lebesgue measurable essentially bounded on $\RR^+$
functions, $g_i(t) \geq 0$, $i=1,\dots, s$,
${\displaystyle 
\int_0^{+\infty} g_i(s)~ds = +\infty}, ~ i=1,\dots,s$.  

\item{{\bf (a5)}} 
$\Phi: (-\infty,0] \to \RR^s$ is a continuous bounded vector function.
\end{description}

Examples of \eqref{intro1} include a system with several concentrated delays
\begin{equation}
\label{ex_syst1}
\frac{dx_i}{dt} = g_i(t) \left[  \sum_{j=1}^{n_i} \alpha_{ij} f_i\left( x_1(h_{i1j}(t)), \dots, x_s(h_{isj}(t))
\right)  - x_i(t) \right], 
\end{equation}
where $h_{ikj}$ satisfy (a2), $\displaystyle \sum_{j=1}^{n_i} \alpha_{ij} = 1$, $i=1, \dots, s$, as well as
a system of integro-differential equations
\begin{equation}
\label{ex_syst2}
\frac{dx_i}{dt} = g_i(t) \left[ 
\sum_{j=1}^{s} \int_{h_{ij}(t)}^t  K_{ij}(t,\tau) f_i(x_1(\tau), \dots, x_s(\tau))~d\tau - x_i(t) \right], 
\end{equation}
with $h_{ij}$ satisfying (a2),
\begin{equation}
\label{ex_syst2_cond}
\sum_{j=1}^{s} \int_{h_{ij}(t)}^t  K_{ij}(t,\tau)~d\tau \equiv 1,~~  K_{ij}(t,\tau) \geq 0, ~~ i,j=1, \dots, s.
\end{equation}

\begin{definition}
\label{def2}
(see \cite{Liz})
Let $D=(a_1,b_1) \times (a_2,b_2) \times \dots \times (a_s,b_s)$, $I_n= [a_{1n},b_{1n}] \times [a_{2n},b_{2n}] \times \dots \times [a_{sn},b_{sn}]$, $F:\RR^s \to \RR^s$.
An equilibrium $z_* \in D$ is {\bf a strong attractor in $D$}  of the difference system
\begin{equation}
\label{dif_syst}
X(n+1)=F(X(n)), ~~n=0,1, \dots 
\end{equation} 
if there exists a sequence of sets $\{I_n\}$, $n=0,1, \dots$,  such that
\begin{equation}
\label{attractor}
{\rm Int}(I_0) = D, ~F(I_n) \subset I_{n+1} \subset {\rm Int}(I_n), ~~n=0,1, \dots, ~~ \bigcap_{n=1}^{+\infty} I_n = z_*.
\end{equation}
\end{definition}

Note that, once $z_*$ is a strong attractor of $F$ in $D$, it is unique, moreover, there are no other equilibrium points of $F$ in $D$.

\begin{uess}
\label{lemma1}
Let $F$ be a continuous function, where $F$ is defined in \eqref{F_def},
and $F(I_k) \subset I_{k+1} \subset {\rm Int}(I_k)$, $F(I_{k+1}) \subset I_{k+2} \subset {\rm Int}(I_{k+1})$,
for some $k \in \NN$, where $I_k= [a_{1k},b_{1k}] \times [a_{2k},b_{2k}] \times \dots \times [a_{sk},b_{sk}]$.

Then there exist $\overline{J}=[\bar{c}_{1},\bar{d}_{1}] \times [\bar{c}_{2},\bar{d}_{2}] 
\times \dots \times [\bar{c}_{s},\bar{d}_{s}]$ and 
$\underline{J}= [{c}_{1},{d}_{1}] \times [{c}_{2},{d}_{2}] \times \dots \times [{c}_{s},{d}_{s}]$ 
such that $\overline{J} \subset {\rm Int}(I_k)$,  $I_{k+1} \subset {\rm Int} \overline{J}$, 
$I_{k+2} \subset {\rm Int} \underline{J}$, $\underline{J} \subset {\rm Int}(I_{k+1})$ and $F(\overline{J}) \subset \underline{J}$.
\end{uess}
{\bf Proof.} For simplicity, we choose $\displaystyle c_i=0.5(a_{i k+1}+a_{i k+2})$, $\displaystyle d_i=0.5(b_{i k+1}+b_{i k+2})$,
then $I_{k+2} \subset$Int$\underline{J}$ and $\underline{J} \subset$Int$(I_{k+1})$ are obviously satisfied.

Next, introduce a family $\{ J_{\alpha} \}$ of compact subsets of the interior of $I_k$ as 
$$
J_{\alpha}= [ \alpha a_{1k} +(1-\alpha) a_{1 k+1}] \times [ \alpha a_{2k} +(1-\alpha) a_{2 k+1}]  \times \dots \times [ \alpha a_{sk} +(1-\alpha) a_{s k+1}], ~~\alpha\in [0,1]
$$
and notice that for $\alpha=0$, $J_{0} = I_{k+1}$, $F(J_{0}) \subset I_{k+2}$, and $I_{k+2} \subset$Int$\underline{J}$. 
Thus there exists $\alpha_0 \in [0,1]$
such that
\begin{equation*}
\alpha_0 = \inf \left\{ \alpha\in [0,1] : F(J_{\alpha}) \subset \underline{J} \right\}, 
\end{equation*}
as the set in the right-hand side is non-empty. If $\alpha_0 >0$, we choose $\alpha=\min\{ \alpha_0, \frac{1}{2} \}$ (to avoid $\alpha=1$) and
denote $\overline{J}=J_{\alpha}$. Then $\overline{J} \subset $Int$(I_k)$,  $I_{k+1} \subset {\rm Int} \overline{J}$, 
$I_{k+2} \subset$Int$\underline{J}$, $\underline{J} \subset $Int$(I_{k+1})$ and $F(\overline{J}) \subset \underline{J}$, and the proof is complete.

It remains to exclude the case $\alpha_0=0$. If for any positive $\alpha$, $F(J_{\alpha}) \not\subset \underline{J}$, we choose a sequence
$\alpha_n=\frac{1}{n}$; by our assumption, there is a sequence of points $z_n \in J_{1/n}$ such that $F(z_n) \not\in \underline{J}$.
By definition, all $z_n\in I_k$ which is a compact set, thus there is a subsequence convergent to some $\bar{z}$, and $F(\bar{z})$ does not belong 
to the interior of \underline{J}. However, as $\alpha_n \to 0$, this limit point $\bar{z}$ belongs to $J_0=I_{k+1}$. However,
$F(I_{k+1}) \subset I_{k+2}$ and Int$(\underline{J}) \subset I_{k+2}$, thus $F(\bar{z})\in $Int$(\underline{J})$, which is a contradiction.
Thus there exists $\alpha_0>0$, and the proof is complete. 
\qed


\section{Main Results}
\label{sec:stability}

Now, we are in a position to prove the main statement of the paper.

\begin{guess}
\label{theorem1}
Let $D=(a_1,b_1) \times (a_2,b_2) \times \dots \times (a_s,b_s)$, $F:\RR^s \to \RR^s$.
Suppose (a1)-(a5) hold, and $z_*$ is a strong attractor of $F$ in $D$, with
$I_n= [a_{1n},b_{1n}] \times [a_{2n},b_{2n}] \times \dots \times [a_{sn},b_{sn}]$, $n\in {\mathbb N}_0 = \{0\}
\cup {\mathbb N}$.

Then for any initial function such that
\begin{equation}
\label{th1eq1}
\Phi \in  C((-\infty, t_0),D), 
\end{equation}
the solution of \eqref{1a},\eqref{2star} satisfies $\lim\limits_{n\to +\infty} X(t)=z_*$.
\end{guess}
{\bf Proof.} We prove that for any initial function satisfying \eqref{th1eq1}, first,
$X(t)\in$Int$(I_0)$ for any $t\geq t_0$ and, in addition, there is $t_1 \geq t_0$ such that $X(t_1) \in$Int$(I_1)$.
Moreover, $X(t) \in$Int$(I_1)$ for any $t \geq t_1$.

Assume that with \eqref{th1eq1} satisfied, there is the first point $t^*$ such that $X(t^*) \not\in \rm{Int}(I_0)$, i.e. $X(t^*)$ is on the boundary $\partial I_0$. Then there exists $i \in \{1, \dots, s\}$ such that either $x_i(t^*)=a_{i0}$ or $x_i(t^*)=b_{i0}$, and $F(X(t))\in I_1$, $f_i(X(t))>a_{i1}$, $t\in [t_0,t^*)$.
In the former case, due to continuity of $x_i$, there is a $t_0^* \in [t_0,t^*)$ such that
$x_i(t) < 0.5( a_{i0} + a_{i1} )$ for $[t_0^*,t^*)$. Hence, using the notation of \eqref{intro1a}, we get
\begin{align*}
x_i(t^*) & = x_i(t_0^*) + \int_{t_0^*}^{t^*} g_i(t)\left[\int_{H(t)} d_{\tau} R_i(t,\tau) f_i(X(\tau)) - x_i(t) \right]~dt 
\\
& >  a_{i0} +  \int_{t_0^*}^{t^*} g_i(t) \left[ a_{i1} \int_{H(t)} d_{\tau} R_i(t,\tau) - \frac{a_{i0} +
a_{i1}}{2} \right] 
\\ & > a_{i0} + 
\int_{t_0^*}^{t^*} g_i(t) \frac{a_{i1} -
a_{i0}}{2} \,dt >a_{i0},
\end{align*}
which contradicts to the assumption $x_i(t^*)=a_{i0}$.

Similarly, in the latter case, assuming that $x_i(t^*)=b_{i0}$ and $x_i(t) > 0.5(b_{i1}+ b_{i0})$
for $t \in [t_1^*,t^*]$, we obtain
\begin{align*}
x_i(t^*) & = x_i(t_1^*) + \int_{t_1^*}^{t^*} g_i(t)\left[\int_{H(t)} d_{\tau} R_i(t,\tau) f_i(X(\tau)) - x_i(t) \right]~dt 
\\ & < b_{i0} +  \int_{t_1^*}^{t^*} g_i(t) \left[ b_{i1} \int_{H(t)} d_{\tau} R_i(t,\tau) - \frac{b_{i1}+b_{i0}}{2} \right] \\ & < b_{i0} + 
\int_{t_1^*}^{t^*} g_i(t) \frac{ b_{i0} - b_{i1} }{2} \, dt < b_{i0},
\end{align*}
again leading to a contradiction. Thus, $X(t)\in$Int$(I_0)$ for any $t\geq t_0$.

Note that by Lemma~\ref{lemma1}, there exist $\overline{J}=[\bar{c}_{1},\bar{d}_{1}] \times [\bar{c}_{2},\bar{d}_{2}] \times \dots \times [\bar{c}_{s},\bar{d}_{s}]$ and $\underline{J}= [{c}_{1},{d}_{1}] \times [{c}_{2},{d}_{2}] \times \dots \times [{c}_{s},{d}_{s}]$ such that $\overline{J} \subset $Int$(I_0)$, $I_1 \subset$Int$\overline{J}$, $I_{2} \subset$Int$\underline{J}$, $\underline{J} \subset $Int$(I_{1})$ and $F(\overline{J}) \subset \underline{J}$.

The proof that there exists $t_1> t_0$, such that $X(t) \in$Int$(I_1)$ for $t \geq t_1$ will consist of two parts. First, we prove that there is a $\bar{t}$ such that $X(\bar{t}) \in$Int$(\overline{J})$; moreover, $X(t)\in $Int$(\overline{J})$ for any $t \geq \bar{t}$. Second, we find $t_1 \geq \bar{t}$ for which $X(t_1) \in$Int$(I_1)$, 
and also justify that $X(t) \in$Int$(I_1)$, $t \geq t_1$.

Since $I_1 \subset$Int$\overline{J}$, we have $\bar{c}_i < a_{i1} < b_{i1} < \bar{d}_i$.
Denote
\begin{equation}
\label{tendrate}
\delta := \min \left\{ \min_{1\leq i \leq s} \left( a_{i1} - \bar{c}_i \right),  
\min_{1\leq i \leq s} \left( \bar{d}_i - b_{i1} \right) \right\}
>0,
\end{equation}
which is positive as a minimum of $2s$ positive values.

For any $X(t) \in$Int$(I_0)$ we have $F(X(t))\in I_1$.
Assume that this does not hold for some $i$.
Let $t$ be such that the $i$th component of $X(t)$ satisfies $x_i \leq \bar{c}_i$. We prove that
there is a moment of time $t_{1i}$ such that  $x_i(t_{1i})  > \bar{c}_i$.

We recall that
\begin{equation}
a_{i0} < \bar{c}_i < a_{i1} < c_i <a_{i2} < b_{i2} < d_i < b_{i1} <  \bar{d}_i < b_{i0}, ~~i=1,\dots,s.
\label{points}
\end{equation}
As long as $x_i(t) \leq \bar{c}_i$,  we have, by \eqref{tendrate},
\begin{align*}
x_i^{\prime}(t)  & = g_i(t) \left[ \int_{H(t)} d_{\tau} R_i(t,\tau) f_i(X(\tau)) - x_i(t) \right]  
\\ & \geq g_i(t) 
\left[ \int_{H(t)} a_{i1}~d_{\tau} R_i(t,\tau) - \bar{c}_i \right] = g_i(t) (a_{i1} - \bar{c}_i) \geq \delta g_i(t). 
\end{align*}
By (a4), the integral of the positive right-hand side diverges, thus there is $t_{i1}$ such that 
$x_i(t_{i1}) >   \bar{c}_i$.
Since there is $t<t_{i1}$ such that
$x_i(t_{i1})\leq \bar{c}_i < a_{i1}$, without loss of generality we can assume $x_i(t) < a_{i1}$, $t<t_{i1}$.
Next, note that $x_i(t) >   \bar{c}_i$ for any $t \geq t_{i1}$. 
In fact, assuming the contrary that $t^*>t_{i1}$ is the smallest value exceeding $t_{i1}$ at which $x(t^*)=\bar{c}_i$,
we get $x_i(t) \in  [\bar{c}_i,a_{i1}]$ for $t\in [t_{i1},t^*]$. Recall that $f_i(X(t)) \geq a_{i1}$, $t\in [t_{i1},t^*]$,
therefore
\begin{align*}
x_i(t^*) = & x(t_{i1})+ \int_{t_{i1}}^{t^*}  g_i(\zeta) 
\left[ \int_{H(\zeta)} d_{\tau} R_i(\zeta,\tau) f_i(X(\tau)) - x_i(\zeta) \right]~d\zeta
\\
\geq & x(t_{i1})+ \int_{t_{i1}}^{t^*}  g_i(\zeta) \left[ a_{i1} - a_{i1} \right]~d\zeta = x(t_{i1}).
\end{align*}
This contradicts to our assumption that $x(t^*)=\bar{c}_i<x(t_{i1})$. 

Similarly, for $x_i(t) \geq \bar{d}_i$,
\begin{align*}
x_i^{\prime}(t) & = g_i(t) \left[ \int_{H(t)}  d_{\tau} R_i(t,\tau) f_i(X(\tau)) - x_i(t) \right]  \\
& \leq g_i(t) 
\left[ \int_{H(t)} b_{i1}~d_{\tau} R_i(t,\tau) - \bar{d}_i \right] =g_i(t) (b_{i1} - \bar{d}_i) \leq - \delta g_i(t). 
\end{align*}
Hence there is $t_{i2}$ such that $x_i(t_{i2}) <   \bar{d}_i$. Again, we justify that $x_i(t) <   \bar{d}_i$ for any $t\geq t_{i2}$.
Choosing
$$
\bar{t}= \max\left\{ \max_{1 \leq i \leq s} t_{i1}, \max_{1 \leq i \leq s} t_{i2} \right\},
$$
we conclude that $X(t)\in $Int$(\overline{J})$ for any $t \geq \bar{t}$.

By (a2), there is $\bar{t}_1$ such that $h_i(t) \geq  \bar{t}$ for any $t \geq \bar{t}_1$, $i=1, \dots, s$.


Note that $F(\overline{J}) \subset \underline{J}$ and $\underline{J} \subset $Int$(I_{1})$.
Define
\begin{equation}
\label{tendrate1}
\delta_1 := \min \left\{ \min_{1\leq i \leq s} \left( c_i - a_{i1}  \right),  \min_{1\leq i \leq s} \left( b_{i1} - d_i \right) \right\},
\end{equation}
which is positive by \eqref{points}.
We have for $t \geq \bar{t}_1$,
$$
x_i^{\prime}(t)  \geq g_i(t) 
\left[ \int_{H(t)} c_i ~d_{\tau} R_i(t,\tau) - a_{i1}  \right] =g_i(t) (c_i - a_{i1}) \geq \delta_1 g_i(t)
$$
for any $x_i(t) \leq a_{i1}$ and conclude by (a4) that there is $t_{i3}$ such that $x_i(t_{i3}) \in (a_{i1}, c_i)$.
Moreover, $x_i(t) >  a_{i1}$ for any $t \geq t_{i3}$. Assuming the contrary that 
$x_i(t) \in  (a_{i1}, x_i(t_{i3})]$ for $t\in [t_{i3}, t_2^*)$ and $x_i(t_2^*)=a_{i1}$, we get for $t\in [t_{i3}, t_2^*)$, 
$f_i(X(t)) \geq c_i$
$$
x_i^{\prime}(t)  \geq g_i(t) 
\left[ \int_{H(t)} c_i ~d_{\tau} R_i(t,\tau) - x_i(t_{i3})  \right] =g_i(t) (c_i - x_i(t_{i3})) 
\geq g_i(t) (c_i - a_i) \geq \delta_1 g_i(t) >0,
$$
which contradicts to the assumption that $x_i(t_2^*)=a_{i1}<x_i(t_{i3})$.

If $x_i(t) \geq b_{i1}$, 
$$
x_i^{\prime}(t) \leq g_i(t) 
\left[ \int_{H(t)} d_i~d_{\tau} R_i(t,\tau) -  b_{i1} \right] =r(t) (d_i - b_{i1}) \leq - \delta g_i(t). 
$$
Thus $x_i (t_{i4}) < b_{i1}$ for some $t_{i4}$,
and similarly we get $x_i(t) <  b_{i1}$ for any $t \geq t_{i4}$. 
Then, for $t \geq t_1$, where
$$
t_1= \max\left\{ \max_{1 \leq i \leq s} t_{i3}, \max_{1 \leq i \leq s} t_{i4} \right\},
$$
we have $x_i(t)\in (a_{i1},b_{i1})$, $i=1,\dots,s$, or $X(t) \in$Int$(I_1)$, $t \geq t_1$.

Since $x_i(t) \in (a_{i1},b_{i1})$, $i=1,\dots, s$, $t \geq t_1$, we proceed to the next induction step from $I_1$ to $I_2$. 

By (a2), there exists $t_1^h > t_0$ such that $h(t)>t_1$ for any $t>t_1^h$. Then, we have an initial value problem with all initial values in Int$(I_1)$ and complete the induction step similarly, justifying that there is $t_2>t_1$ such that $X(t)\in I_2$ for $t \geq t_2$, and also $h(t)>t_2$ for any $t>t_2^*$. Proceeding in the same manner from $n$ to $n+1$, we prove that there is an increasing sequence of $t_n$ such that $X(t) \in I_n$ for $t \geq t_n$. Since the intersection of $I_n$ is $z_*$,
this implies $\displaystyle \lim_{t\to +\infty} X(t)=z_*$.
\qed

Consider system \eqref{intro1} or \eqref{1a} under the assumptions (a1)-(a5), as well as an additional assumption
\begin{description}
\item{{\bf (a6)}}
$z_*=(x_1^*, \dots, x_j^*, \dots x_s^*)$ is the only equilibrium in the domain $D$,
and there exist $L_{ij} \geq 0$ such that,
for any $(x_1, \dots, x_s) \in D$, a.e.
\begin{equation}
\label{L_coefficients}
\begin{array}{ll}
\displaystyle
& \left| f_i(x_1, \dots, x_{j-1}, x_j,x_{j+1}, \dots, x_s)- f_i(x_1, \dots, x_{j-1}, x_j^*,x_{j+1}, \dots, x_s)\right|  \vspace{2mm }  \\
\leq &  L_{ij} \left| x_j - x_j^* \right|,~~  i,j=1, \dots, s, ~~ L = \left(  L_{ij} \right)_{i,j=1}^s.
\end{array}
\end{equation}
\end{description}

Note that \eqref{L_coefficients} in (a6) is satisfied if $f_i$ is globally Lipschitz, i.e.
for any $(x_1, \dots,x_j,\dots, x_s) \in D$ and $ (x_1, \dots,y_j,\dots, x_s)\in D$,  a.e.,
$$\displaystyle
\left| f_i(x_1, \dots, x_{j-1}, x_j,x_{j+1}, \dots, x_s)- f_i(x_1, \dots, x_{j-1}, y_j,x_{j+1}, \dots, x_s)\right|
\leq L_{ij} \left| x_j - y_j \right| $$
for $i,j=1, \dots, s$. 
In particular, if $f_i$ are a.e. differentiable and
\begin{equation*}
\left| \frac{\partial f_i}{\partial x_j} \right| \leq L_{ij} \mbox{ a.e. },
~~i,j=1, \dots, s,
\end{equation*}
condition (a6) is satisfied.

We recall that a matrix $A=(a_{ij})_{i,j=1}^s$ is {\em nonnegative} if $a_{ij} \geq 0$ and {\em positive} if $a_{ij}>0$, $i,j=1, \dots,s$.
Let $\| X \|$ be an arbitrary fixed norm of a column vector in $\RR^s$, and $\| A \|$ be the induced matrix norm.
The classical definition of an $M-$matrix will be used. Following \cite{Berman}, we say that $A=(a_{ij})_{i,j=1}^s$ is a (non-singular) {\em $M$-matrix} if $a_{ij}\leq 0$ for $i\neq j$ and $A^{-1}$ is positive.
By $I$ we denote an $s \times s$ identity matrix.


There are many equivalent definitions of $M$-matrices, see \cite{Berman} and also \cite[Lemma 2.3]{BAM2014}.

\begin{uess}\cite[p.~137,142, Exercise 2.9 of Chapter 6]{Berman}\label{lemma_Berman}
$A$ is an $M-$matrix if and only if $a_{ij}\leq 0$, $i\neq j$ and 
there exist positive numbers $\xi_i$, $i=1,\dots, s$ such that $$\xi_i a_{ii}>\sum_{j\neq
i}\xi_j|a_{ij}|,~ i=1,\dots, s.$$
\end{uess}

\begin{uess}
\label{lem_M_matrix}
Let (a6) be satisfied, $L$ be defined in \eqref{L_coefficients}, and $I-L$ be an $M$-matrix.
Then, there exist $\{a_{in}\}_{n=1}^{+\infty}$ and  $\{b_{in}\}_{n=1}^{+\infty}$, $i=1, \dots,s$ such that 
$$a_{in} < a_{i n+1}< x_i^* < b_{i n+1} < b _{in}, ~~ n\in {\mathbb N}, ~~i=1, \dots, s,$$
$$\lim_{n\to +\infty} a_{in}= \lim_{n\to +\infty} b_{in} = x_i^*, ~~i=1, \dots, s,$$
and the domains 
$$I_n=[a_{1n}, b_{1n}] \times [a_{2n}, b_{2n}] \dots \times [a_{sn}, b_{sn}]$$
satisfy
$$
I_1 \subset D, ~~ F(I_n) \subset I_{n+1}, ~~n \in {\mathbb N}.
$$
\end{uess}
{\bf Proof.} 
By definition of the nonnegative matrix $L$, all off-diagonal entries of $I-L$ are non-positive.
By Lemma~\ref{lemma_Berman}, we have a finite set of $\xi_i>0$, $i=1,\dots, s$ such that 
$$
\xi_i \left( 1- L_{ii} \right) > \sum_{j\neq i} \xi_j L_{ij},~~i=1, \dots, s,
$$
or
$$
\alpha_i := \sum_{j\neq i} \frac{\xi_j}{\xi_i} L_{ij} + L_{ii} <1 ,~~i=1, \dots, s.
$$
Denote 
\begin{equation}
\label{alpha}
\alpha= \max_{1 \leq i \leq s} \alpha_i \in (0,1). 
\end{equation}
Choose for some $c >0$,
$$a_{i1}=x_{i}^* -  c \xi_i,~~b_{i1} = x_{i}^* + c \xi_i,~~ i=1, \dots, s$$
such that
$$
 I_1=[a_{11}, b_{11}] \times [a_{21}, b_{21}] \dots \times [a_{s1}, b_{s1}] \subset D.
$$
In particular, if 
$$
D= (a_{10}, b_{10}) \times (a_{20}, b_{20}) \times \dots \times (a_{s0}, b_{s0}),
$$
we can take any positive $c$ satisfying
$$
c \leq \min \left\{ \min_{1 \leq j \leq s} \frac{b_{j0}-x_j^*}{\xi_j}  , \min_{1 \leq j \leq s} \frac{x_j^* - a_{j0}}{\xi_j}\right\}.
$$
We have
\begin{equation}
\label{d1_def}
I_1=[x_1^*-c \xi_1, x_1^*+c \xi_1] \times [x_2^*-c \xi_2, x_2^*+c \xi_2] \times \dots \times [x_s^*-c \xi_s, x_s^*+c \xi_s]
\end{equation}
implying $|x_j-x_j^*|< c\xi_j$, $j=1, \dots, s$ for $X\in I_1$ and $f_i(x_1^*,x_2^*, \dots, x_s^*)= x_i^*$.
By \eqref{L_coefficients}, once $X\in I_1$,
\begin{align*}
\left| f_i(X) - x_i^* \right|  \leq & \left| f_i(x_1,x_2, \dots, x_{s-1}, x_s)-f_i(x_1^*,x_2, \dots,x_{s-1}, x_s) \right|
\\ & +
\left| f_i(x_1^*,x_2, \dots, x_{s-1},x_s) - f_i(x_1^*,x_2^*, \dots, x_{s-1},x_s) \right| \\ &  + \dots 
+ \left| f_i(x_1^*,x_2^*, \dots, x_{s-1}^*,x_s) 
- f_i(x_1^*,x_2^*, \dots, x_{s-1}^*,x_s^*) \right| \\
\leq & L_{ii} \left| x_i-x_i^* \right| + \sum_{j\neq i}  L_{ij} \left| x_j-x_j^*  \right| \\
\leq &  L_{ii} c \xi_i + \sum_{j\neq i}  L_{ij} c \xi_j \leq \alpha c \xi_i,
\end{align*}
where $\alpha \in (0,1)$ is denoted in \eqref{alpha}.

Recall \eqref{d1_def} and denote for $n \in {\mathbb N}$,
\begin{equation}
\label{dn_def}
I_{n+1}=[x_1^*- \alpha^{n} c \xi_1, x_1^*+ \alpha^{n} c \xi_1] \times 
\dots \times [x_s^*-\alpha^{n} c \xi_s, x_s^*+ \alpha^{n} c \xi_s].
\end{equation}

We have justified $F(I_1) \subset I_2$, with $I_2$ defined in \eqref{dn_def}. 
Now let $X \in I_n$. Then,
\begin{align*}
\left| f_i(X) - x_i^* \right|  \leq & \left| f_i(x_1,x_2, \dots, x_{s-1},x_s)-f_i(x_1^*,x_2, \dots, x_{s-1}, x_s) \right| \\ & +
\left| f_i(x_1^*,x_2, \dots, x_{s-1},x_s) - f_i(x_1^*,x_2^*, \dots, x_{s-1},x_s) \right| \\ &  + \dots 
+ \left| f_i(x_1^*,x_2^*, \dots, x_{s-1}^*, x_s) - f_i(x_1^*,x_2^*, \dots, x_{s-1}^*,x_s^*) \right| \\
\leq & L_{ii} \left| x_i-x_i^* \right| + \sum_{j\neq i}  L_{ij} \left| x_j-x_j^*  \right| \\
\leq & L_{ii} \alpha^{n-1} c \xi_i + \sum_{j\neq i}  L_{ij} \alpha^{n-1} c \xi_j \leq \alpha  \alpha^{n-1}c \xi_i
= \alpha^{n}c \xi_i,
\end{align*}
so $F(X)\in I_{n+1}$, where $I_{n+1}$ is defined in \eqref{dn_def}.
$F$ is an $\alpha$-contraction, $\alpha \in (0,1)$. 
Since
$$
a_{in}= x_{i}^* -  \alpha^{n-1} c \xi_i,~~b_{in} = x_{i}^* + \alpha^{n-1} c \xi_i, ~~i=1, \dots, s,
$$
we get $\lim\limits_{n\to +\infty} a_{in}= \lim\limits_{n\to +\infty} b_{in} = x_i^*$, $i=1, \dots, s$,
which concludes the proof of the lemma.
\qed

Lemma~\ref{lem_M_matrix} and Theorem~\ref{theorem1} immediately imply the following asymptotic stability result.

\begin{guess}
\label{prop2}
Suppose (a1)-(a6) are satisfied.
If $I-L$ is an $M$-matrix, where $L$ is defined in \eqref{L_coefficients}, then any solution of \eqref{intro1} with $X_0 \in D$ converges to $z_*$.
\end{guess}
{\bf Proof.} Let $I-L$ be an $M$-matrix. By Lemma~\ref{lem_M_matrix} and Definition~\ref{def2}, $z_*$ is a strong attractor in $D$ of difference system \eqref{dif_syst}. 
Thus, the conditions of Theorem~\ref{theorem1} are satisfied, and therefore any solution of \eqref{intro1},\eqref{2star} with \eqref{th1eq1} being fulfilled 
converges to $z_*$. 
\qed

\section{Applications and Examples}
\label{sec:applications}

Consider a particular case of $s=2$.  System (\ref{1a}) includes the model with variable delays
\begin{equation}
\label{1avar}
\begin{array}{ll}
\displaystyle \frac{dx}{dt} & \displaystyle  = g_1(t) \left[ f_1(y(h_1(t))) - x(t) \right],
\vspace{2mm} \\
\displaystyle \frac{dy}{dt} & \displaystyle  = g_2(t) \left[ f_2(x(h_2(t))) - y(t) \right], ~~t \geq 0
\end{array}
\end{equation}
and the integro-differential system 
\begin{equation}
\label{1aintegro}
\begin{array}{ll}
\displaystyle \frac{dx}{dt} & \displaystyle  = g_1(t) \left[ \int_{h_1(t)}^t K_1 (t,s)f_1(y(s))~ds  - x(t) \right],
\vspace{2mm} \\  
\displaystyle \frac{dy}{dt} & \displaystyle  = g_2(t) \left[ \int_{h_2(t)}^t K_2 (t,s)f_2(x(s))~ds  - y(t) \right],~~t \geq 0,
\end{array}
\end{equation}
where for both \eqref{1avar} and \eqref{1aintegro}, the functions $h_i$ and  $g_i$ satisfy (a2) and (a4), respectively. For \eqref{1aintegro}, in addition, a modification of (a3)
\begin{description}
\item{{\bf (a3$^*$)}}   
$K_i(t, s):\RR^+ \times \RR^+ \to \RR^+$, $i=1,2$ are locally integrable functions
in both $t$ and $s$ satisfying $\displaystyle \int_{h_i(t)}^t K_i (t,s)~ds \equiv 1$, $i=1,2$
\end{description}
is assumed to hold.

Further, for both \eqref{1avar} and \eqref{1aintegro}, the functions $f_i$ should satisfy
\begin{description}
\item{{\bf (a7)}}
Both $f_1:\RR^+ \to \RR^+$ and $f_2:\RR^+ \to \RR^+$ are continuous strictly monotone increasing, $f_1(0)=f_2(0)=0$ and 
\begin{equation}
\label{two_syst_inter}
f_2(x)>f_1^{-1}(x), ~~x \in (0,x^*), ~~~ f_2(x)<f_1^{-1}(x), ~~x \in (x^*,+\infty),~y^*=f_1(x^*).
\end{equation}
\end{description}
Note that \eqref{two_syst_inter} implies that  $(x^*,y^*)$ is the unique equilibrium of systems \eqref{1avar} and \eqref{1aintegro}.

\begin{prop}
\label{prop1}
Let $h_i$, $g_i$ and $f_i$ satisfy (a2),(a4) and (a7), respectively, and, in the case of  \eqref{1aintegro}, (a3$^*$) hold.
Then, all solutions of \eqref{1avar} and \eqref{1aintegro} with non-negative 
non-trivial in both $x$ and $y$ continuous initial conditions converge to $(x^*,y^*)$.
\end{prop}
{\bf Proof.} 
Define $D=(0,+\infty) \times (0,+\infty)$, $F=(f_1,f_2)^T$, then $F:D \to D$.
Let us construct $I_n=[a_{1n},b_{1n}] \times [a_{2n}, b_{2n}]$, where
\begin{equation}
\label{a_sequence}
a_{11}<a_{12} \dots < a_{1n}<a_{1 n+1} < \dots, ~~
a_{21}<a_{22} \dots < a_{2n}<a_{2 n+1} < \dots ,
\end{equation}
\begin{equation}
\label{b_sequence}
\dots <b_{1 n+1} < b_{1n} < \dots b_{12}<b_{11}, ~~
\dots <b_{2 n+1 } < b_{2n} < \dots b_{22}<b_{21},
\end{equation}
and 
\begin{equation}
\label{ab_limits}
\lim_{n\to +\infty} a_{1n} = \lim_{n\to +\infty} b_{1n} = x^*,~~
\lim_{n\to +\infty} a_{2n} = \lim_{n\to +\infty} b_{2n} =y^*,
\end{equation}
which would imply that $(x^*,y^*)$ is a strong attractor. 

Since $f_1$ is monotone increasing, so is $f_1^{-1}$, also both $f_2$ and $f_2^{-1}$ are monotone increasing. By \eqref{two_syst_inter} in (a6),
$f_2(x) > f_1^{-1}(x)$ for $x\in (0,x^*)$. 
Denote $y=f_2(x)$, $x=f_2^{-1}(y)$.
The function $f_1$ is also monotone increasing, thus, taking $f_2$ of both sides, we get
$$ f_1(f_2(x))>f_1(f_1^{-1}(x))=x, \mbox{~~ or ~~~} f_1(y)>f_2^{-1}(y), ~~y \in (0,y^*).$$
Similarly, considering $x>x^*$, or $y=f_2(x)>y^*$, we get $f_1(y)<f_2^{-1}(y)$ for $y \in (y^*,+\infty)$.
Thus
$$
f_1(y)>f_2^{-1}(y), ~~y \in (0,y^*), ~~~ f_2(y)<f_1^{-1}(y), ~~y \in (y^*,+\infty).
$$
Next, choose arbitrary initial left bounds $a_{11}\in (0,x^*)$ and $b_{11} \in (x^*,+\infty)$.
For the left bound define $a_{21}=f_2(a_{11})$, $a_{12}=f_1(a_{21})$,
$a_{22}=f_2(a_{12})$.
By \eqref{two_syst_inter},
$$
0<f_1^{-1}(x)<f_2(x)<f_2(x^*)=y^*, ~~x \in (0,x^*), 
$$
hence $f_2: (0,x^*) \to (0,y^*)$. 
Recall that $f_2$ is monotone and $a_{11}\in (0,x^*)$, therefore
$a_{21}=f_1(a_{11}) \in (0,y^*)$. 
In addition, for $x\in  (0,x^*)$,  \eqref{two_syst_inter} implies $f_1: (0,y^*) \to (0,x^*)$ for monotone increasing $f_1$.
Therefore $a_{12} \in (0,x^*)$ and $a_{22}=f_2(a_{12}) \in (0,y^*)$. We have
$a_{21}, a_{22} \in (0,y^*)$, $a_{12} \in (0,x^*)$. 
Also,
$$
a_{12}=f_1(a_{21}) = f_1(f_2(a_{11})) > f_1(f_1^{-1}(a_{11}))=a_{11},
$$
$$
a_{22}=f_2(a_{12}) = f_2(f_1(a_{21})) > f_2(f_2^{-1}(a_{21})) =a_{21}.
$$
For an induction step, take
\begin{equation}
\label{ind_step}
a_{1 n+1} = f_1(a_{2n}), \quad a_{2 n+1} = f_2(a_{1 n+1}).
\end{equation}
From $a_{1n} \in (0,x^*)$, $a_{2n} \in (0,y^*)$ and monotonicity of $f_1$,$f_2$ we get
$a_{1 n+1} \in (0,x^*)$, $a_{2 n+1} \in (0,y^*)$, as well as
$$
a_{1 n+1}=f_1(a_{2n}) = f_1(f_2(a_{1n})) > f_1(f_1^{-1}(a_{1n}))=a_{1n},
$$
$$
a_{22}=f_2(a_{1 n+1}) = f_2(f_1(a_{2n})) > f_2(f_2^{-1}(a_{2n})) =a_{2n}.
$$
Thus, \eqref{a_sequence} holds, and we have two monotone increasing sequences $\{ a_{1n}\}$ and $\{ a_{2n}\}$ 
bounded by $x^*$ and $y^*$, respectively, from above.
Hence both sequences have limits $\displaystyle \lim_{n\to +\infty} a_{1n} = d_1 \in (0,x^*]$, $\displaystyle \lim_{n\to+\infty} a_{2n} = d_2 \in (0,y^*]$. By \eqref{ind_step} and continuity of $f_1,f_2$, $d_1=f_1(d_2)$, $d_2=f_2(d_1)$, which implies $d_1=x^*$, $d_2=y^*$. 

For the right bound we use 
$$
y^*=f_2(x^*) < f_2(x)<f_1^{-1}(x), ~~x \in (x^*,+\infty),
$$
$f_2:(x^*,+\infty) \to (y^*,+\infty)$ and $f_2:(y^*,+\infty) \to (x^*,+\infty)$.
Therefore, we get bounds for
$b_{21}=f_2(b_{11}) \in (y^*,+\infty)$, $b_{12}=f_1(a_{21}) \in (x^*,+\infty)$,
$b_{22}=f_1(b_{12}) \in (y^*,+\infty)$.

The sequences of $b_{1n}\in (x^*,+\infty)$, $b_{2n} \in (y^*,+\infty)$ satisfying \eqref{b_sequence} are constructed similarly
$$
b_{1 n+1} = f_1(b_{2n}), \quad b_{2 n+1} = f_2(b_{1 n+1}),
$$
and the proof of \eqref{b_sequence} follows the same steps, as well as $\displaystyle \lim_{n\to +\infty} b_{1n} = x^*$, $\displaystyle \lim_{n\to +\infty} b_{2n} =y^*$.
Therefore, \eqref{ab_limits} is satisfied, and $(x^*,y^*)$ is a strong attractor. Thus, all the conditions of Theorem~\ref{theorem1} are satisfied.
The application of Theorem~\ref{theorem1} concludes the proof. 
\qed

The statement of Proposition~\ref{prop1} is the main result of \cite{Nonlin2015},
and a two-dimensional cooperative system described in \cite{Nonlin2015} is a particular case of the system considered in the present paper.


\begin{example}
Consider a particular case of \eqref{ex_syst1} 
\begin{equation}
\label{ex_syst3}
\begin{array}{lll}
x' & = & g_1(t) \left[ \sqrt{y(h_1(t))} - x(t) \right], \vspace{2mm} \\
y' & = & g_2(t) \left[ \sqrt{x(h_2(t))} - y(t) \right],~~t \geq 0.
\end{array}
\end{equation}
For \eqref{ex_syst3}, $f_1(y)=\sqrt{y}$, $f_1^{-1}(x)=x^2$, $f_2(x)=\sqrt{x}$,
$f_2^{-1}(y)=y^2$, $f_1(1)=1$, $f_2(1)=1$  and 
$$
f_2(x) =\sqrt{x} > f_1^{-1}(x)=x^2,~~x\in (0,1), ~~ f_2(x) =\sqrt{x} < f_1^{-1}(x)=x^2,~~x\in (1,+\infty).
$$
Thus \eqref{two_syst_inter} holds, and
Proposition~\ref{prop1} implies that any solution with nonnegative nontrivial initial conditions converges to (1,1).
\end{example}

\begin{example}
For the system
\begin{equation}
\label{ex_syst6}
\begin{array}{lll}
x' & = & \displaystyle g_1(t) \left[ y^2(h_1(t)) - x(t) \right], \vspace{2mm} \\
y' & = & \displaystyle g_2(t) \left[ \sqrt[4]{x(h_2(t))} - y(t) \right],~~ t \geq 0,
\end{array}
\end{equation}
the functions $f_1(y)=y^2$ and $f_2(x)=\sqrt[4]{x}$ are continuous and monotone increasing on ${\mathbb R}_+$, $f_1(1)=1$, $f_2(1)=1$. Also, 
$$
f_2(x) =\sqrt[4]{x} > f_1^{-1}(x)=\sqrt{x},~~x\in (0,1), ~~ f_2(x) =\sqrt[4]{x} < f_1^{-1}(x)=\sqrt{x},~~x\in (1,+\infty).
$$
Since \eqref{two_syst_inter} holds, by Proposition~\ref{prop1} any non-negative non-trivial in both $x$ and $y$ solution converges to $(1,1)$.
\end{example}

\begin{example}
\label{BAM}

For the BAM neural network without delays in the leakage terms 
\begin{equation*}
x_i^{\prime} = g_i(t) \left[ \left( \sum_{j=1}^s \alpha_{ij} x_j(h_{ij}(t)) \right)^{1/(2k_i)} - x_i(t) \right],
\end{equation*}
where
\begin{equation*}
\alpha_{ij}  \geq 0, ~~ \sum_{j=1}^s \alpha_{ij} =1, ~~k_i \in \NN,~~i,j =1, \dots s, ~~t \geq 0,
\end{equation*}
Theorem~\ref{prop2}  implies that the equilibrium $(1,1, \dots, 1)$ attracts all solutions with non-negative non-trivial continuous initial conditions.
\end{example}

\begin{remark}
In \cite{BAM2014}, a neural system which can be reduced to 
$$
\dot{x}_i(t) = \alpha_i(t) \left[ -x_i(h_i(t))+ \sum_{j=1}^s F_{ij}(t,x_j(h_{ij}(t)) \right],~~t \geq 0, ~~i=1,\dots,s
$$
was considered. If in the leakage terms $h_i(t) \equiv t$, $i=1,\dots,s$, the results of  \cite[Theorem 2.5]{BAM2014} coincide with a particular case of Theorem~\ref{prop2} when delays are concentrated. Thus,
compared to \cite[Theorem 2.5]{BAM2014},
Theorem~\ref{prop2} considers more general distributed delays in non-leakage part but assumes a particular case of
non-delayed  leakage terms. Therefore, the results are independent.
\end{remark}

\begin{remark}
According to \cite[Theorem 3.1]{Liz}, all solutions of the system
$$
\dot{x}_i(t) = -x_i(t)+ \sum_{j=1}^s \alpha_{ij}(t)f_j(x_j(t-\tau_{ij})),~~t>0, ~~i=1,\dots,s,
$$
where
\begin{equation}
\label{flower}
\left| f_j(u)-f_j(v) \right| \leq L_j |u-v|, ~\forall u,v \in {\mathbb R},~~ \max_{1\leq i \leq s} \sum_{j=1}^s \alpha_{ij} L_j <1,
\end{equation}
converge to the zero equilibrium $(0,0,\dots,0)$.
Following the notation of the present paper, denote by $L$ the matrix with the entries $L_{ij}=\alpha_{ij} L_j$.
Thus zero is globally attractive, once a sum of the entries of each column is less than one.
Note that Theorem~\ref{prop2} states attractivity of the zero equilibrium once the matrix $I-L$ is an $M$-matrix.
For example, let
$$ L = \left[ \begin{array}{cc} \frac{1}{2} & 2 \vspace{2mm} \\ \frac{1}{16} & \frac{1}{2} \end{array} \right],$$
then \eqref{flower} is not satisfied,
since the sum of the entries of the second column exceeds 1, but it is easy to check that 
$I-L$ is an $M$-matrix by its form and positivity of the inverse matrix
$$I-L=\left[ \begin{array}{cc} \frac{1}{2} & -2 \vspace{2mm} \\ -\frac{1}{16} & \frac{1}{2} \end{array} \right],
~~(I-L)^{-1}= \left[ \begin{array}{cc} 4 & 16\vspace{2mm}  \\ \frac{1}{2} & 4 \end{array} \right],
$$
therefore Theorem~\ref{prop2} implies global attractivity of the zero equilibrium.
\end{remark}

Next, consider the Nicholson-type system  
\begin{equation}
\label{Nicholson}
\frac{dx_i}{dt} = g_i \left[\sum_{j \neq i} a_{ij} x_j +
\sum_{k=1}^m \beta_{ik} x_i(\tau_{ik}(t)) e^{-x_i(\tau_{ik}(t))} - x_i(t)\right],~~ i=1, \dots, s,~~t \geq 0,
\end{equation}
where $g_i>0$, $a_{ij}$ and $\beta_{ik}$ are non-negative, while
$$
\beta_i := \sum_{k=1}^m \beta_{ik} >0,
$$
and for some $\tau>0$, ~ $t-\tau_{ik}(t) \leq \tau,~~i=1,\dots, s,~k=1, \dots, m$.
Global attractivity conditions for \eqref{Nicholson} were obtained in \cite{Faria2,Faria1}, 
see also references therein and  \cite{Faria2} for a detailed history outline. A positive equilibrium for this system exists \cite{Faria2} once all the constants
\begin{equation}
\label{gamma}
\gamma_i := \frac{\beta_i}{1 - \sum_{j\neq i} a_{ij}}
\end{equation}
satisfy $\gamma_i >1$, $i=1, \dots, s$.

The unique positive equilibrium of \eqref{Nicholson} exists and is globally attractive 
\cite[Theorems 2.5 and 3.3]{Faria2}, once
\begin{equation}
\label{abs_Nichol}
1< \gamma_i < e^2, ~~i=1, \dots,s.
\end{equation}
Note that inequalities \eqref{abs_Nichol} imply $1-\sum_{j\neq i} a_{ij}< \beta_i < e^2$, $i=1, \dots,s$.

As an application of our results,
consider the system 
\begin{equation}
\label{1}
\frac{dx_i}{dt} = g_i(t) \left[\sum_{j \neq i} a_{ij} x_j(t) + 
\int\limits_{h_i(t)}^{t} \beta_{i} x_i(\tau) e^{-x_i(\tau)}~d_{\tau} r_i(t,\tau) - x_i(t)\right], ~i=1, \dots, s,~~t \geq 0,
\end{equation}
where for the functions $g_i, r_i$ conditions (a2)-(a4) hold, $a_{ij}\geq 0, 1<\beta_{i}\leq e^2$.

In particular, \eqref{1}
includes the system with several concentrated delays generalizing \eqref{Nicholson}
\begin{equation*}
\frac{dx_i}{dt} = g_i(t) \left[\sum_{j \neq i} a_{ij} x_j(t) + 
\sum_{k=1}^{m_i} \beta_{ik} x_i(\tau_{ik}(t)) e^{-x_i(\tau_{ik}(t))} - x_i(t)\right], ~~i=1, \dots, s,~~t \geq 0,
\end{equation*}
where $\displaystyle \beta_i = \sum_{k=1}^{m_i} \beta_{ik}>0$, $i=1, \dots, s$, $g_i$ satisfy (a4), for $\tau_{ik}$ condition (a2) holds.

Assume that $x^*=(x_1^*,\dots,x_s^*)$ is a unique positive equilibrium of \eqref{1}.
In particular, the fact that $\gamma_i>1$, where $\gamma_i$ are defined in \eqref{gamma}, $i=1, \dots, s$  guarantees that such an equilibrium exists, similarly to systems with concentrated delays. 

Denote 
\begin{equation}
\label{star1}
\alpha_i=\left\{\begin{array}{ll}
\max\{1-\ln \beta_i, \beta_i e^{-2}\},& 1<\beta_i\leq e,\\
\beta_i e^{-2},& e<\beta_i<e^2. 
\end{array}\right.
\end{equation}

\begin{guess}
\label{theorem_Nicholson}
Let (a2)-(a4) hold and $I-L$ be an $M$-matrix, where 
\begin{equation*}
L=(L_{ij})_{i,j=1}^s,~L_{ij}=\left\{\begin{array}{ll}
a_{ij},& j\neq i,\\
\alpha_i,& j=i,\end{array}\right.
\end{equation*}
and $\alpha_i$ are denoted in \eqref{star1}.
Then all solutions of \eqref{1} with non-negative non-trivial initial conditions converge to $x^*$.
\end{guess}
{\bf Proof.} Denote 
$$
f_i(x_1, \dots, x_s)=\sum_{j \neq i} a_{ij} x_j + \beta_{i} x_i e^{-x_i}.
$$
To apply Theorem~\ref{prop2}, we have to estimate the partial derivatives 
$\displaystyle \left|\frac{\partial f_i}{\partial x_j}\right|$.
We have 
$$
\left|\frac{\partial f_i}{\partial x_j}\right|=\left\{\begin{array}{ll}
a_{ij},& j\neq i,\\
\beta_i|1-x_i|e^{-x_i},& j=i.\end{array}\right.
$$

The maximum of the function $xe^{-x}$ is attained at $x=1$ and equals $1/e$.
According to \cite[Theorem 2.6]{BrZhuk}, any positive solution of the equation
\begin{equation}
\label{2}
\frac{dy_i}{dt} = g_i(t) \left[\int_{h_i(t)}^{t}  \beta_{i} y_i(\tau) e^{-y_i(\tau)} ~d_{\tau} r_i(t,\tau) -  y_i(t)\right]
\end{equation}
with $\beta_i>1$ satisfying 
$$
\limsup_{t\to +\infty} y_i(t) \leq \sup_{x\in[0, +\infty)} \beta_i xe^{-x}=\frac{\beta_i}{e}, $$
$$
\liminf_{t\to +\infty} y_i(t) \geq \min_{x\in (1,+\infty)} \beta_i xe^{-x} = x_i^0,
$$
where
$$x_i^0:= \left\{ \begin{array}{ll}  \ln \beta_i, &  \beta_i \in [1,e), \\  \frac{\beta_i^2}{e}e^{-\frac{\beta_i}{e}},
&  \beta_i \in [e,e^2). \end{array} \right. $$

Next, let $x_i$ be a component of a solution in \eqref{1}. Then, with the same initial conditions as in \eqref{2}, 
since all components are positive, $ x_i(t)\geq x_i^0$.
Note that a similar result for concentrated delays was justified in  \cite[Theorem 2.3]{BIT}.
Hence it is sufficient to estimate $p_i(x_i):=\beta_i|1-x_i|e^{-x_i}$ only on the interval $[x_i^0,+\infty)$.
There are two cases: $x_i^0  \leq 1$ corresponding to $\beta_i \in (1,e]$ and $x_i^0 >1$  for $\beta_i \in (e,e^2)$.

If $x_i^0 \leq x_i<1$ then $p_i=\beta_i(1-x_i)e^{-x_i}$, $p_i^{'}=-\beta_i(2-x_i)e^{-x_i}<0$.
Hence  $\max\limits_{x_i^0 \leq x_i<1} p_i(x_i)=p_i(x_i^0)=\beta_i(1-x_i^0)e^{-x_i^0} = 1-\ln \beta_i$. 

If $x_i^0<1$, $x_i>1$ then $p_i=\beta_i(x_i-1)e^{-x_i}, p_i^{'}=\beta_i(2-x_i)e^{-x_i}$
and $\max\limits_{x_i>1} p_i(x_i)=p_i(2)=\beta_ie^{-2}$.

If $x_i \geq x_i^0>1$ then $\max_{x_i \geq x_i^0>1} p_i(x_i)=p_i(2)=\beta_ie^{-2}$.

Theorem~\ref{prop2} implies that, since the matrix $I-L$ is an $M$-matrix,  $x^*$
is a global attractor, which concludes the proof.
\qed


In particular, if $\beta_i>e$, the diagonal entries of $I-L$ are positive, and the matrix is diagonally dominant
\begin{equation*}
\beta_i e^{-2} < 1- \sum_{i \neq j} a_{ij}, ~~i=1, \dots, s,
\end{equation*}
the unique positive equilibrium of \eqref{1} is globally asymptotically stable, which generalizes the right inequality in \eqref{abs_Nichol} to the case of distributed delays and variable growth rates.

\begin{corol}
\label{corol_Nichol}
Let $n=2$, (a2)-(a4) hold  and 
\begin{equation}
\label{2cond}
e< \beta_1 < e^2,~e<\beta_2 < e^2, ~~a_{12}a_{21}<(1-\alpha_1)(1-\alpha_2),
\end{equation}
where $\alpha_i$ are introduced in \eqref{star1}.
Then the positive equilibrium is globally attractive.
\end{corol}
{\bf Proof.} 
For $n=2$, as $\alpha_i=\beta_i e^{-2}$, we have
$$
I-A=\left(\begin{array}{ll}
1-\alpha_1&-a_{12}\\
-a_{21}&1-\alpha_2.
\end{array}\right) = \left(\begin{array}{ll}
1-\beta_1 e^{-2} &-a_{12}\\
-a_{21}&1-\beta_2 e^{-2}
\end{array}\right) .
$$
Thus $I-A$ is an $M$-matrix if $\beta_i e^{-2} <1$, $i=1,2$ and $a_{12}a_{21}<(1-\alpha_1)(1-\alpha_2)$, in particular, when 
\eqref{2cond} holds.
\qed


\begin{example}
\label{example_Nichol_2}
Consider the system with $h_1$, $h_2$ satisfying (a2), $r_1>0$, $r_2>0$,
\begin{equation}
\label{system_2_Nichol}
\begin{array}{lll}
x'(t) & = & \displaystyle r_1 \left[ 0.5 y(t)+ 4x(h_1(t))e^{-x(h_1(t))} -x(t)  \right], \\
y'(t) & = & \displaystyle r_2 \left[ 0.2 x(t)+ 5y(h_2(t))e^{-y(h_2(t))} -y(t)  \right] ,~~t \geq 0.
\end{array}
\end{equation}
Obviously $\beta_1=4$ and $\beta_2=5$ are in $(e,e^2)$. Also, 
$$a_{12}a_{21} =0.1 <(1-\beta_1 e^{-2})(1-\beta_2 e^{-2}) \approx 0.148295,$$
so \eqref{2cond} is satisfied, and the positive equilibrium is globally attractive.

Note that for $n=2$, conditions \eqref{abs_Nichol} are equivalent to
\begin{equation}
\label{abs_Nichol_2}
1-a_{12} < \beta_1 < (1-a_{12}) e^2, ~~ 1-a_{21} < \beta_2 < (1-a_{21}) e^2.
\end{equation}
The right inequalities in \eqref{abs_Nichol_2} can be rewritten as 
$$
a_{12}<1-\beta_1 e^{-2}, ~~ a_{21}< 1- \beta_2 e^{-2},
$$
where the first inequality is not satisfied since
$$a_{12}=0.5 > 1-\beta_1 e^{-2} \approx 0.45866.$$
Thus Corollary~\ref{corol_Nichol} establishes global attractivity of the positive equilibrium of \eqref{system_2_Nichol}, while \eqref{abs_Nichol_2} fails.
\end{example}

\section{Discussion}
\label{sec:discussion}

General system \eqref{1a} was motivated by neural networks but another common application is a compartment, or patch model of mathematical biology. For example, \eqref{1} is a particular type of a compartment model, where $x_i$ is  a population size in the $i$th patch,
$a_{ij}(t)$ describes the relocation rate from the patch $j$ to patch $i$, $i \neq j$, and Nicholson's growth rate. Assuming the logistic growth rate, we get for $K_i>0$ being the carrying capacity of the $i$th patch, a model
\begin{equation}
\label{1abcd}
\frac{dx_i}{dt} = g_i(t) \left[\sum_{j \neq i} a_{ij} x_j(t) + 
\int\limits_{h_i(t)}^{t} \beta_{i} x_i(\tau) \left( 1 - \frac{x_i(\tau)}{K_i} \right)~d_{\tau} r_i(t,\tau) - x_i(t)\right], ~i=1, \dots, s.
\end{equation}
As possible extension of current research, another compartment model with the Mackey-Glass growth rate
\label{Mackey2}
\begin{equation}
\label{Mack2}
\frac{dx_i}{dt} = r_i(t) \left[\sum_{j \neq i} a_{ij} x_j + 
\int\limits_{h_i(t)}^{t} \beta_{i} \frac{x_i(\tau)} {1+x_i^n(\tau)} ~d_{\tau} r_i(t,\tau) -x_i(t)\right], i=1, \dots, s
\end{equation}
can be explored under usual assumptions. 
It would be interesting to investigate existence, uniqueness and  absolute attractivity of the positive equilibrium, and the dependency of this equilibrium on the parameters, as well as delay-dependent stability.

In addition to Nicholson-type system \eqref{1} studied in the present paper and 
proposed \eqref{1abcd}, \eqref{Mack2}, it is possible to consider Ricker-type model, for $i=1, \dots, s$,
\begin{equation}
\label{Ricker}
\frac{dx_i}{dt} = g_i(t) \left[ 
\int\limits_{h_i(t)}^{t} \beta_{i} x_i(\tau) \exp \left\{K_i-x_i(\tau) - \sum_{j \neq i} a_{ij} x_j(\tau) \right\}~d_{\tau} r_i(t,\tau) - x_i(t)\right]. 
\end{equation}
Global attractivity of a positive equilibrium for $s=2$ and $s=3$ was recently studied in \cite{Cabral}, with explicit criteria obtained. It would be interesting to compare sufficient conditions under which the positive equilibrium of \eqref{Ricker} attracts all positive solutions with these tests. In general, the strong attractivity is a stricter assumption that the fact that all solutions of a system of difference equations converge to a certain solution \cite{Liz2}, so it is expected that global attractivity conditions for \eqref{Ricker} may be more restrictive than the tests in \cite{Cabral}.


Let us discuss whether we can replace a sequence of parallelepiped-type domains containing a fixed point $z_*$ by any closed compact sets
including $z_*$.
Notice that a compact set on a line mentioned in the definition of a strong attractor in \cite{Liz} can be reduced to a closed segment such that its interior contains $z_*$. Recall that every open set $\RR$ is a union of at most countable number of open disjoint intervals \cite[p. 45, problems 22 and 29]{Rudin}. Hence a closed bounded set is a union of at most countable number of disjoint closed segments (some may consist of one point only). As segments are disjoint, only one of the segments includes $z_*$. Therefore at each stage we can consider only this segment. The fact that this segment has a non-empty interior, follows from \eqref{attractor}.
Thus, instead of a sequence of compact sets in \cite{Liz}, without loss of generality we can consider 
$I_n= [a_{1n},b_{1n}] \times 
\dots \times [a_{sn},b_{sn}]$ as in Definition~\ref{def2}. 
Thus, our definition in fact coincides with the one in \cite{Liz}.

The main result of the present paper is the proof of global attractivity of non-autonomous equations with distributed and finite, not necessarily bounded, delays. One of the natural questions arising will be extension of the present results to equations with infinite, but exponentially decaying memory. Considering delay-dependent attractivity conditions for systems with distributed delays, similarly to the case of ``small delays'' in \cite{BrZhuk,Faria2}, is another important question, once a cooperative system is not globally asymptotically stable for any delays. 

The results of the present paper are concerned with non-autonomous systems. For relevant autonomous equations with distributed delays, it has been recently proved \cite{belair,Bernard} that, once we replace a distributed delay in an autonomous equation with its expected value, and the resulting delay equation is stable, so is the model with a distributed delay. It is an interesting and challenging problem to extend this result to autonomous systems with distributed delays.    

\ack
The second author was partially supported by the NSERC grant RGPIN/05976-2015.


\bigskip


\begin{thebibliography}{99}

\bibitem{Aliseyko}
A.\,N. Aliseyko, 
Lyapunov matrices for neutral time-delay systems with exponential kernel, {\em Systems Control Lett.} {\bf 131} (2019), 104497, 7 pp.


\bibitem{Cabral}
E.\,C. Balreira, S. Elaydi, R. Luís, 
Global stability of higher dimensional monotone maps, {\em J. Difference Equ. Appl.} {\bf 23} (2017), 
2037–-2071. 


\bibitem{belair}
S. Bernard, J. B\'{e}lair, M.C. Mackey, Sufficient conditions for
stability of linear differential equations with distributed delay, {\em Discrete Contin.
Dyn. Syst. Ser.B} {\bf 1} (2001),
233-256.

\bibitem{Automatica2012}
L. Berezansky and E. Braverman,
On nonoscillation and stability for systems of differential equations with a 
distributed delay, {\em Automatica J. IFAC} {\bf 48} (2012), 
612-–618. 

\bibitem{Nonlin2013}
L. Berezansky and E. Braverman, 
Stability of equations with a distributed delay,
monotone production and nonlinear mortality, {\em Nonlinearity} {\bf 26} (2013), 
2833–-2849.


\bibitem{JMAA2014}
L. Berezansky and E. Braverman,
On multistability of equations with a distributed delay, monotone production
and the Allee effect,
{\em J. Math. Anal. Appl.} {\bf 415} (2014),
873--888.


\bibitem{Nonlin2015}
L. Berezansky and E. Braverman, 
On stability of cooperative and hereditary systems with a distributed delay, 
{\em Nonlinearity} {\bf 28} (2015), 
1745--1760.

\bibitem{ZAA2019}
L. Berezansky and E. Braverman, 
On stability of delay equations with positive and negative coefficients with applications, 
{\em Z. Anal. Anwend.} {\bf 38} (2019), 
157--189.

\bibitem{BAM2014}
L. Berezansky, E. Braverman and L. Idels,
New global exponential stability criteria for nonlinear delay 
differential systems with applications to BAM neural networks,
{\em Appl. Math. Comput.} {\bf 243} (2014),
899--910.

\bibitem{BIT}
L. Berezansky, L. Idels and L. Troib, 
Global dynamics of Nicholson-type delay systems with applications, {\em Nonlinear Anal. Real World Appl.} {\bf 12} (2011), 
436-–445. 


\bibitem{Berman}
A. Berman and R. Plemmons,  Nonnegative Matrices in the Mathematical Sciences, Classics in Applied Mathematics,  
SIAM, Philadelphia, 1994.

\bibitem{Bernard} 
S.~Bernard and F. Crauste,  
Optimal linear stability condition for scalar differential equations with distributed delay, 
{\em Discrete  Contin. Dyn. Syst. B} {\bf 20} (2015), 1855--1876.

\bibitem{BrZhuk}
E. Braverman and S. Zhukovskiy, 
Absolute and delay-dependent stability of equations with a distributed delay,
{\em Discrete Contin. Dyn. Syst.} {\bf 32} (2012), 
2041-–2061.

\bibitem{Faria2}
D. Caetano and T. Faria, 
Stability and attractivity for Nicholson systems with time-dependent delays, 
{\em Electron. J. Qual. Theory Differ. Equ.} 2017, Paper No. 63, 19 pp.

\bibitem{Campbell}
S.\,A. Campbell and I. Ncube, 
Stability in a scalar differential equation with multiple, distributed time delays,
{\em J. Math. Anal. Appl.} {\bf 450} (2017), 
1104--1122.


\bibitem{Cheng}
Z. Cheng, Y. Wang and J. Cao, 
Stability and Hopf bifurcation of a neural network model with distributed delays and strong kernel,
{\em Nonlinear Dynam.} {\bf 86} (2016), 
323-–335.

\bibitem{Oliveira}
S. Esteves, E. G\"{o}kmen and J.\,J. Oliveira, Global exponential stability of nonautonomous neural network models 
with continuous distributed delays, {\em Appl. Math. Comput.} {\bf 219} (2013), 
9296--9307.

\bibitem{Faria_JDE_2017}
T. Faria, 
Periodic solutions for a non-monotone family of delayed differential equations with applications to Nicholson systems,
{\em J. Differential Equations} {\bf 263} (2017), 
509--533. 

\bibitem{Faria1}
T. Faria and G. R\"{o}st, 
Persistence, permanence and global stability for an n-dimensional Nicholson system, 
{\em J. Dynam. Differential Equations} {\bf 26} (2014), 
723--744. 

\bibitem{Glizer}
V.\,Y. Glizer, Uniform stabilizability of parameter-dependent systems with state and control delays by smooth-gain controls, {\em J. Optim. Theory Appl.} {\bf 183} (2019), 
50--65.

\bibitem{Gourley}
S.\,A. Gourley, R. Liu and Y. Lou, 
Intra-specific competition and insect larval development: a model with time-dependent delay, 
{\em Proc. Roy. Soc. Edinburgh Sect. A} {\bf 147} (2017), 
353--369. 

\bibitem{Hattaf}
K. Hattaf and N. Yousfi, 
A class of delayed viral infection models with general incidence rate and adaptive immune response. Int. J. Dyn. Control {\bf 4} (2016), no. 3, 254--265.

\bibitem{Hopfield}
J.J. Hopfield, Neural networks with graded response have collective computation 
properties like those of two-state neurons, {\em Proc. Natl. Acad. Sci.}
{\bf 81} (1984), 3088--3092.

\bibitem{Beretta}
S. Liu and E. Beretta, Competitive systems with stage structure of distributed-delay type, {\em J. Math. Anal. 
Appl.} {\bf 323} (2006), 
331--343. 

\bibitem{Liu}
X. Liu and P. Stechlinski, Hybrid control of impulsive systems with distributed delays, {\em Nonlinear Anal. Hybrid 
Syst.} {\bf 11} (2014), 57--70. 

\bibitem{Liz}
E. Liz and A. Ruiz-Herrera, Attractivity, multistability, and bifurcation in delayed Hopfield's model with 
non-monotonic feedback, {\em J. Differential Equations} {\bf 255} (2013), no. 11, 4244--4266. 


\bibitem{Liz2}
E. Liz and A. Ruiz-Herrera, Alfonso Addendum to "Attractivity, multistability, and bifurcation in delayed Hopfield's model with non-monotonic feedback'', 
{\em J. Differential Equations} {\bf  257} (2014), 
1307–-1309.
 

\bibitem{Muroya}
Y. Muroya and T. Faria, 
Attractivity of saturated equilibria for Lotka-Volterra systems with infinite delays and feedback controls,
{\em Discrete Contin. Dyn. Syst. Ser. B} {\bf 24} (2019), 
3089-–3114.


\bibitem{Rudin}
W. Rudin, 
Principles of Mathematical Analysis, third edition, 
International Series in Pure and Applied Mathematics, McGraw-Hill Book Co., 
New York-Auckland-Düsseldorf, 1976. 

\bibitem{Solomon}
O. Solomon and E. Fridman, New stability conditions for systems with distributed delays, {\em  Automatica J. IFAC}
{\bf 49} (2013), 
3467--3475. 

\bibitem{Xu}
X. Xu, L. Liu and G. Feng, 
Semi-global stabilization of linear systems with distributed infinite input delays and actuator saturations. Automatica J. IFAC {\bf 107} (2019), 398--405.

\bibitem{Yang}
Y. Yang, L. Zou and S. Ruan, 
Global dynamics of a delayed within-host viral infection model with both virus-to-cell and cell-to-cell transmissions. 
Math. Biosci. {\bf 270} (2015), part B, 183--191. 

\bibitem{Yuan_2011}
Y. Yuan and J. B\'{e}lair, Stability and Hopf bifurcation analysis for functional 
differential equation with distributed delay, {\em SIAM J. Appl. Dyn. Syst.} {\bf  10} (2011), 
551--581. 

\end{thebibliography}
\end{document}